\documentclass[11pt]{article}
\usepackage{amsfonts}
\usepackage{graphicx}        % standard LaTeX graphics tool
\usepackage{amsmath,amssymb,mathrsfs}
\usepackage{amsbsy,amsthm,amsfonts}
\usepackage{amsthm}
\usepackage{color}
\usepackage{fullpage}
\usepackage{cite}
\usepackage[colorlinks,citecolor=blue]{hyperref}

\newcommand\figref[1]{Fig.~\ref{#1}}

\newcommand\sectref[1]{Section~\ref{#1}}

\newcommand{\Omegarm}   {\mathrm{\Omega}}

\newcommand{\bfE}   {\mathbf{E}}
\newcommand{\bfD}   {\mathbf{D}}

\newcommand{\etal}  {\emph{et~al.}}
\newcommand{\bfn}   {\mathbf{n}}

\newcommand{\bfr}   {\mathbf{r}}

\newcommand{\bfP}   {\mathbf{P}}

\newcommand{\calL}    {\mathcal{L}}

\newcommand{\calE}    {\mathcal{E}}

\newcommand{\epsloc}   {\epsilon_{l}} % {\epsilon_{\mathrm{loc}}}
\newcommand{\epsnl}   {\epsilon_{\mathrm{nl}}}
\newcommand{\Omeganl}   {\Omegarm_{\mathrm{nl}}}
\newcommand{\Omegaloc}   {\Omegarm_{l}}

\newtheorem{theorem}{Theorem}

\graphicspath{{Figures/}}

\begin{document}

\begin{center}
\textbf{\Large Trefftz Functions for Nonlocal Electrostatics}
		\vskip 4mm

		Igor Tsukerman

Department of Electrical and Computer Engineering,\\
		The University of Akron, OH 44325-3904, USA\\
		igor@uakron.edu

\today
\end{center}
	
\begin{abstract}
Electrostatic interactions in solvents play a major role 
in biophysical systems. There is a consensus in the literature 
that the dielectric response of aqueous solutions is nonlocal:
polarization depends on the electric field not only at a given point 
but in the vicinity of that point as well. This is typically modeled
via a convolution of the electric field with an appropriate
integral kernel.

A primary problem with nonlocal models is high computational 
cost. A secondary problem is restriction of convolution integrals to the 
solvent, as opposed to their evaluation over the whole space. 

The paper develops a computational tool 
alleviating the ``curse of nonlocality'' and helping to handle the 
integration correctly. This tool is 
\textit{Trefftz approximations}, which tend to furnish much higher 
accuracy than traditional polynomial ones. In the paper, Trefftz approximations 
are developed for problems of nonlocal electrostatics, with the goal
of numerically ``localizing''  the original nonlocal problem.
This approach can be extended to nonlocal problems in other areas of 
computational mathematics, physics and engineering.
\end{abstract}
	
%\maketitle

%%%%%%%%%%%%%%%%%%%%%%%%%%%%%%%%%%%%%%%%%%%%%%%%%
\section{Introduction}
\label{sec:Objectives}
%%%%%%%%%%%%%%%%%%%%%%%%%%%%%%%%%%%%%%%%%%%%%%%%%
%
Electrostatic interactions are well known to play a major role in biomolecular
and biophysical systems. For example, Ren \textit{et al}. write in their review 
paper \cite{Ren-RevBiophys12}, 
\begin{quotation}
	``Among the various components of molecular interactions, 
	[electrostatic interactions]
	are of special importance because of their long-range nature and
	their influence on polar or charged molecules, including water, aqueous 
	ions, proteins,	nucleic acids, carbohydrates, and membrane lipids. 
	In particular, robust models of	electrostatic interactions are essential 
	to understand the solvation properties of biomolecules
	and the effects of solvation upon biomolecular folding, binding, enzyme 
	catalysis, and dynamics.''
\end{quotation}

Central in these models is an accurate representation of the dielectric
properties of aqueous solutions. As e.g. Bardhan \textit{et al}. note 
\cite{Bardhan15},
\begin{quotation}
``	One of the long-standing challenges in molecular biophysics is the 
	development of accurate, yet simple models
	for the influence of biological fluids (aqueous solutions composed of water 
	and dissolved ions) on biological
	molecules such as proteins and DNA.''
\end{quotation}
There is a consensus in the literature that accurate models
of water and aqueous solvents must account for \textit{nonlocality}.
That is, in contrast to local descriptions, polarization $\bfP(\bfr)$ and
the displacement vector field $\bfD(\bfr)$ depend on the
electric field $\bfE$ not only at a given point $\bfr$, but
in the vicinity of that point as well.
This dependence is usually expressed as convolution with an appropriate 
kernel (\sectref{sec:Nonlocal Electrostatics}). 
The significance of nonlocal effects in the electrostatic response of 
solvents has been emphasized since the pioneering work of 
Dogonadze, Kornyshev, Vorotyntsev \textit{et al.} in the 1970s 
\cite{Dogonadze74,Kornyshev80,Bopp96}; see 
\cite{Hildebrandt04,Weggler-Rutka-Hildebrandt-JCompPhys10,
Rubinstein-Sherman-Biopolymers07,Xie-SIAM12,Xie-SIAM13,Bardhan-JChemPhys11,Bardhan12}.

A primary problem with nonlocal models is high computational cost. 
In the Finite element (FE) or Finite Difference (FD) context, system matrices 
become less sparse
by a factor of $\mathcal{O}(\delta / h)^d$, where $\delta$ is the
scale of nonlocal interactions, $h$ is the mesh size, and $d$ is the number
of spatial dimensions. Note that dependence of the computational cost
on the number of nonzero matrix entries is almost always superlinear.

A secondary problem is that
convolution with the $\bfE$ field should be confined to the nonlocal region
(solvent); extending this convolution to the solute with local 
characteristics or to the whole space, as done in many
existing publications, may lead
to qualitatively incorrect results (\sectref{sec:Domain-of-Convolution}).

The approach outlined in the paper is aimed at solving both principal problems:
    removing the computational ``curse of nonlocality''
    and handling the convolution integrals properly.
The main tool for that is \textit{Trefftz approximations}, 
which in many instances furnish much higher accuracy than traditional 
piecewise-polynomial ones. This development is inspired by our successful 
application of Trefftz approximations in a variety of disparate problems: 
wave propagation and scattering, colloidal particles, photonic band structure, 
homogenization of periodic heterostructures 
\cite{Tsukerman-book20,Tsukerman-CAMWA19,Dai-Webb11,Mansha-OpEx17}.
At the same time, construction of Trefftz funtions for nonlocal problems is new 
(\sectref{sec:Trefftz-funcs-nonlocal}).

By definition, Trefftz functions satisfy locally (in weak form) the underlying 
differential equations of the problem. Examples of such functions
are harmonic polynomials for the Laplace equation; plane, cylindrical
or spherical waves for the Helmholtz or Maxwell equations, etc.
The superior accuracy of Trefftz approximations 
explains their growing popularity in a large variety of methods and 
applications: 
Domain Decomposition \cite{Herrera00,Farhat-DD-DG09},
Generalized FEM \cite{Melenk96,Babuska97,Babuska-GFEM2004,Plaks03,Proekt02,
	Strouboulis-GFEM-Helmholtz2006}, 
Discontinuous Galerkin (DG)
\cite{Farhat-DD-DG09,Cockburn00,Arnold02,Buffa-Monk-ultraweak08,Gittelson09,
	Gabard-wave-based-DG-UW-LS11,Hiptmair-PWDG-Helmholtz2011,Kretzschmar-IMA16},
and FD (``Flexible Local Approximation MEthods,'' FLAME) 
\cite{Tsukerman-JCP10,Tsukerman05,Tsukerman06,Tsukerman-PBG08,Tsukerman-JOPA09}.

Trefftz functions derived in this paper can be incorporated into any of
the above methods; the work on nonlocal Trefftz-FLAME and 
Trefftz-DG is planned.
The key idea of FLAME is to replace the Taylor expansions on which
classical FD is based with Trefftz approximations. This often leads
to high-order schemes even in the presence of material interfaces
\cite{Tsukerman05,Tsukerman06,Tsukerman-CAMWA19,
	Tsukerman-PBG08,Tsukerman-JCP10,Tsukerman-PLA17,Tsukerman-book20}.

DG-Trefftz methods possess many attractive features.
Most relevant to the subject of this paper are exponential accuracy
in many cases (\sectref{sec:Harmonic-polynomials}, 
\cite{Hiptmair-Trefftz-survey16,Kretzschmar14,Kretzschmar-IMA16});
adaptivity~\cite{Hiptmair-PWDG-Helmholtz2011,Hiptmair-PWDG16};
a natural treatment of discontinuities -- which is instrumental, in 
particular, for modeling solute-solvent interfaces.
In general, DG methods can handle complex geometries, 
curved boundaries and various boundary conditions.
Equations of \textit{Trefftz}-DG reduce to the skeleton of the mesh,
making the computation more efficient, since volume integrals 
vanish~\cite{Kretzschmar14,Kretzschmar-IMA16,Egger-space-time-DG15}.
%
%%%%%%%%%%%%%%%%%%%%%%%%%%%%%%%%%%%%%%%%%%%%%%%%%%%%%%%%%%%%%%%%%%%%%%%%%%
\section{Nonlocal Electrostatics: Physical Models and Computational Complexity}
\label{sec:Nonlocal Electrostatics}
%%%%%%%%%%%%%%%%%%%%%%%%%%%%%%%%%%%%%%%%%%%%%%%%%%%%%%%%%%%%%%%%%%%%%%%%%%
%
%%%%%%%%%%%%%%%%%%%%%%%%%%%%%%%%%%%%%
\subsection{Physical and Mathematical Models}
\label{sec:Physical-Models}
%%%%%%%%%%%%%%%%%%%%%%%%%%%%%%%%%%%%%
%
\subsubsection{Nonlocal Dielectric Properties}
Nonlocal models, whereby one field at a given point depends on another field
\textit{in the vicinity} of that point, have been introduced in
many areas of physics: plasticity, mechanical vibrations,
liquid crystals, optics, nanostructures
\cite{Jirasek98,Engelen-plasticity03,Shen-nonlocal-plate10,Abdelouhab-nonlocal-waves89,
Ansari-nonlocal-vibrations-graphene10,Jirasek-Rolshoven03,
Park-Calderer-nonlocal-LC08,McMahon09,Agranovich84}.
In the latter area, nonlocal treatment is particularly common in the 
description of
plasmonic effects and in homogenization of metamaterials, where nonlocality
often goes under the name of ``spatial dispersion'' (due to
dependence of material parameters on the wave vector in Fourier space).

In macromolecular simulaiton, the importance of taking nonlocality 
into account is now widely recognized.
As noted in the introduction, nonlocal electrostatic models were developed 
by Dogonadze, Kornyshev, Vorotyntsev and collaborators 
in the 1970s 
\cite{Dogonadze74,Kornyshev80,Bopp96} and since then have been widely
used and extended
\cite{Basilevsky-Parsons96,Hildebrandt04,Weggler-Rutka-Hildebrandt-JCompPhys10,
	Rubinstein-Sherman-Biophys04,Ren-RevBiophys12,Rubinstein-Sherman-Biopolymers07,
	Xie-SIAM12,Xie-SIAM13,Bardhan-JChemPhys11,David-Abajo11,Bardhan12}.

In linear \textit{local} models the displacement field $\bfD$ is related to the
electric field $\bfE$ as
\begin{equation}\label{eqn:D-eq-eps-E}
	\bfD(\bfr) \,=\, \epsloc (\bfr) \bfE(\bfr)
\end{equation}
where $\epsloc$ is the dielectric permittivity of the medium and $\bfr$
is the position vector; linearity is assmued throughout the proposal.
For simplicity, we also assume isotropy, which is accurate for aqueous solutions
in the bulk, although the analysis and methods
can be extended to anisotropic media if necessary.

It is important to distinguish convolution over the whole space
$\mathbb{R}^d$ ($d = 1,2,3$) and the convolution-like integration
\textit{restricted to a given domain} $\Omegarm$; we denote
the former and the latter with `*' and  `$*_{\Omegarm}$', respectively:
\begin{equation}\label{eqn:u-conv-Omega-v-def}
	\calE *_{\Omegarm} \bfE ~~\overset{\mathrm{def}}{=}~
	\int_{\Omegarm}  \calE(\bfr - \bfr') \, \bfE(\bfr') d \bfr'
\end{equation}
The difference between integration over a bounded domain
and the whole space does play a significant role
(\sectref{sec:Domain-of-Convolution}).

Nonlocal constitutive relation is assumed to have the form\footnote{Although 
the local term could 
	be absorbed into the convolution integral as a Dirac delta function,
	it is expedient not to do so.}
\begin{equation}\label{eqn:D-eq-eps-E-plus-Eps-E}
	\bfD(\bfr) \,=\, \epsloc (\bfr) \bfE(\bfr) \,+\,
	\epsnl \, \calE *_{\Omegarm_{\mathrm{nl}}} \bfE
\end{equation}
where $\Omegarm_{\mathrm{nl}}$ is the domain occupied by the nonlocal
medium, and $\epsnl$ is a given parameter. 
%
%%%%%%%%%%%%%%%%%%%%%%%%%%%%%%%%%%%%%
\subsubsection{The Nonlocal Poisson Equation}
\label{sec:Nonlocal-Poisson}
%%%%%%%%%%%%%%%%%%%%%%%%%%%%%%%%%%%%%
%
In the chosen simply connected computational domain $\Omegarm$
with a Lipschitz-continuous boundary, 
the electrostatic potential $u$, the electric and displacement fields 
$\bfE$ and $\bfD$ satisfy
\begin{equation}\label{eqn:ES-equations}
	\bfE = -\nabla u; \quad \nabla \cdot \bfD = \gamma \rho;
	\quad 
	u \in H^1(\Omegarm)
\end{equation}
where $\rho$ is a given charge density and the coefficient $\gamma$ depends on 
the system of units.

Let $\Omegarm$ be partitioned into subdomains $\Omegarm_l$ and $\Omeganl$
occupied by local and nonlocal media, where the respective constitutive 
relation are \eqref{eqn:D-eq-eps-E} and \eqref{eqn:D-eq-eps-E-plus-Eps-E}.
%
%%
%\begin{equation}\label{eqn:D-eq-Eps-E}
%	\bfD(\bfr) = \epsloc(\bfr) \bfE(\bfr), ~~~ \bfr \in \Omegarm_l;
%	\quad
%	\bfD(\bfr) = \calE *_{\Omegarm_{nl}} \bfE 
%	\overset{\mathrm{def}}{=} \int_{\Omegarm_{nl}}
%	\calE(\bfr - \bfr') \bfE(\bfr') \, d\bfr'
%\end{equation}
%
On the interface boundary between $\Omegarm_l$ and $\Omeganl$,
the standard continuity condition for the normal component of $\bfD$ holds:
\begin{equation}\label{eqn:Dn-continuous}
	\bfD_{\mathrm{nl}} \cdot \bfn_{\mathrm{nl}}
	+ 	\bfD_l \cdot \bfn_l = 0
\end{equation}
where $\bfn_{\mathrm{nl}}$ and $\bfn_l$ are the outward unit normals
to $\Omegarm_l$ and $\Omeganl$, respectively, and the rest of the notation
is self-explanatory.
%
%%%%%%%%%%%%%%%%%%%%%%%%%%%%%%%%%%%%%
\subsubsection{The Poisson-Boltzmann Equation}
\label{sec:PBE}
%%%%%%%%%%%%%%%%%%%%%%%%%%%%%%%%%%%%%
%
Consider an electrolyte with positive and negative salt ions 
carrying equal and opposite charges. In general, there may be 
several species of ions; their statistical distribution
in an electrostatic potential $u$ leads to the Poisson--Boltzmann equation
(PBE) 
\cite{Gouy10,Chapman13,Fogolari01,Fogolari02,Deserno:011401:2002}:\footnote{Unit-dependent
	coefficients such as $4\pi$ or $\epsilon_0$ are for brevity omitted.}
\begin{equation}\label{eqn:Poisson-Boltzmann}
	\epsilon_{\mathrm{s}} \nabla^2 u ~=~ -\rho ~-~ \sum_{\alpha} n_{\alpha} 
	q_{\alpha}
	\exp \left( -\frac{q_{\alpha} u}{k_B T} \right)
\end{equation}
where summation is over all species of ions present in the solvent,
$n_{\alpha}$ is volume concentration of species $\alpha$ in the
bulk, $q_{\alpha} = Z_\alpha e$ is the charge of species $\alpha$;
$k_B$ is the Boltzmann constant, $T$ is temperature, and $\rho$
is the density of all charges other than the microions.
If the electrostatic energy $q_{\alpha} u$ is much smaller than the
thermal energy $k_B T$, then the PBE can be approximately linearized around
$u = 0$ to yield
\begin{equation}\label{eqn:Poisson-Boltzmann-linearized-simple}
	\epsilon_{\mathrm{s}} \nabla^2 u \,-\, ( k_B T)^{-1}
	\left( \sum_{\alpha} n_{\alpha} q_{\alpha}^2 \right) u ~=~ 0,
\end{equation}
%
%($\kappa$ is called the \emph{Debye--H\"{u}ckel parameter}).
%
%It is useful to estimate the order of magnitude of the potential for
%which linearization is acceptable. Equating electrostatic energy $q
%u$ of monovalent ions ($q = e$) to thermal energy $k_B T$, one
%obtains the threshold $u_{kT} = k_B T / e$ $\approx$~25~mV at room
%temperature.
%
%Equation \eqref{eqn:Poisson-Boltzmann-linearized-simple} is known
%as the Debye--H\"{u}ckel approximation. 
The potential satisfying this equation will typically exhibit exponential 
decay away from the sources.
%with the characteristic length (the \emph{Debye--H\"{u}ckel 
%length}).
%equal to the inverse of $\kappa$.
%
Eq.~\eqref{eqn:Poisson-Boltzmann-linearized-simple} can be extended
to the nonlocal case by including the convolution term
$\epsnl \nabla \cdot (\calE *_{\Omegarm_{\mathrm{nl}}} \nabla u$). 
%
%%%%%%%%%%%%%%%%%%%%%%%%%%%%
\subsubsection{Limitations of the PBE Model}\label{sect:Limitations-PBE}
%%%%%%%%%%%%%%%%%%%%%%%%%%%%
%
The main physical assumption behind the PBE is that each mobile
charge is effectively in the \emph{mean field} of all other charges,
and has the Boltzmann probability of acquiring any given energy.
This probability is assumed to be \emph{unconditional}, i.e. not
depending on possible redistribution of other ions in response to
the motion of a given ion. In other words, mean field theory
disregards any correlations between the positions and movement of 
the ions. However, it is demonstrated in \cite{Nguyen00}
that such correlations may in fact be appreciable for multivalent ions.
A consensus exists that at least for monovalent ions
the correlations are weak enough for the PB model to be valid,
and the linearized PBE \eqref{eqn:Poisson-Boltzmann-linearized-simple}
then provides a physically accurate description.
The approach of this paper is extendable to the nonlocal version
of the linearized PBE; Trefftz functions for that equation
can be derived following the same procedure as for the nonlocal Poisson equation
(\sectref{sec:Trefftz-funcs-nonlocal}).
%
%%%%%%%%%%%%%%%%%%%%%%%%%%%%
\subsubsection{Solvent--Solute Boundary Value 
Problems (Nonlocal)}\label{sect:Solvent-solute-bvp}
%%%%%%%%%%%%%%%%%%%%%%%%%%%%
%
The formulation of solvent--solute boundary value problems is standard.
Inside the solute (such as a macromolecule), the electrostatic potential
is assumed to satisfy the classical Poisson equation 
$\epsilon_m \nabla^2 u(\bfr) = -\gamma \rho_m(\bfr)$, where $\epsilon_m$
is the (relative) dielectric permittivity of the molecule 
(typically assumed to be in the range from $\sim$2 to 4), $\rho_m(\bfr)$
is the atomic charge density of the molecule (assumed to be 
given),\footnote{The electronic charge density is accounted for via 
polarization -- that is, via $\epsilon_m$.} 
and $\gamma = 4\pi$ in the Gaussian system or 
$\epsilon_0^{-1}$ in SI.

Within the solvent (domain $\Omegarm_{\mathrm{nl}}$), 
in which the solute is immersed, and in the absence of charges, 
the electrostatic potential satisfies the nonlocal Laplace equation 
\begin{equation}\label{eqn:nonlocal-Laplace}
	\nabla \cdot (\calE *_{\Omegarm_{\mathrm{nl}}} \nabla u) 
	\,\equiv\, 
	\nabla \cdot \int_{\Omegarm_{\mathrm{nl}}}
	\calE(\bfr - \bfr') \, \nabla u(\bfr') \, d \bfr'
	~=~ 0
\end{equation}
or, alternatively, the nonlocal Poisson-Boltzmann equation if microions are 
present.
Across the solute-solvent interface, the potential is continuous,
and so is the normal component of the $\bfD$ field. These boundary
conditions are conventional, but in the solvent $\bfD$ is related
to $\nabla u$ nonlocally.

Numerical solution of the solvent--solute problem is facilitated by
potential splitting, as $\rho_m$ contains point charges which
are difficult to handle numerically otherwise. Such splitting
of the solution into the ``forced'' and ``reaction'' parts --
in this case, $u = u_{\rho} + u_r$ -- is standard in the theory
of differential equations, in wave scattering problems 
(incident + scattered fields), in magnetostatics \cite[Chap 
10]{Tsukerman-book20}, and in many other instances; $ u_{\rho}$
is the Coulomb potential of the given charges. Rather than applying
the potential splitting in the whole space, it may be conceptually
and algorithmically simpler to confine it to the solute.
The resultant boundary conditions on the solute-solvent interface can be easily 
handled, especially in DG.
%
%%%%%%%%%%%%%%%%%%%%%%%%%%%%%%%%%%%%%
\subsection{The Domain of Convolution Makes a Significant Difference}
\label{sec:Domain-of-Convolution}
%%%%%%%%%%%%%%%%%%%%%%%%%%%%%%%%%%%%%
%
In most numerical procedures for the nonlocal electrostatic problem to date,
the integration for the $\bfD$ field in \eqref{eqn:nonlocal-Laplace} is 
extended from the solvent domain $\Omegarm_s$ to the whole space, 
with the tacit physical assumption that this does not significantly
affect the results. This extension is not just a matter of practical
convenience; rather, it plays a critical role in the simplification
of the problem.

Indeed, convolution-like integration over a fixed bounded domain,
in contrast with integration over the whole space,
does not possess several key properties and hence does not
easily lend itself to simplification. First, ``finite-domain convolution''
is not in general commutative, a trivial example of which is
\begin{equation}\label{eqn:1D-conv-1-star-x}
	1 *_{[0,1]} x \,\equiv \, \int_0^1 1 \cdot x' dx'
	= \frac12 
	~~~\neq~~~
	x *_{[0,1]} 1 \,\equiv \, \int_0^1 (x-x') \cdot 1 dx'= x - \frac12
\end{equation}
Secondly, differentiation of  finite-domain convolution does not
in general reduce to differentiating either of the terms; e.g.
\begin{equation}\label{eqn:diff_1D-conv}
	d_x (	x *_{[0,1]} 1) \,=\, d_x \left( x - \frac12 \right) = 1
	~~~\neq~~~ x *_{[0,1]} d_x 1 = 0	
\end{equation}
Yet the localization method, proposed originally by
Hildebrandt \cite{Hildebrandt04,Weggler-Rutka-Hildebrandt-JCompPhys10} and later
developed, enhanced and implemented by others 
\cite{Xie-SIAM12,Xie-SIAM13,Ying-Xie18},
does rely on the standard differentiation rules for full-space
convolutions (\sectref{sec:Hildebrandt-Localization}).

\vspace{1mm}
\textit{Remark}. The bounded region of integration could be rigorously
expanded to the whole space if the solution is padded with zero
outside the computational domain. In that case, however,
differentiation will be polluted by distributional derivatives
(boundary delta functions). 
\vspace{1mm}

\begin{figure}
	\centering
	\includegraphics[width=0.48\linewidth]{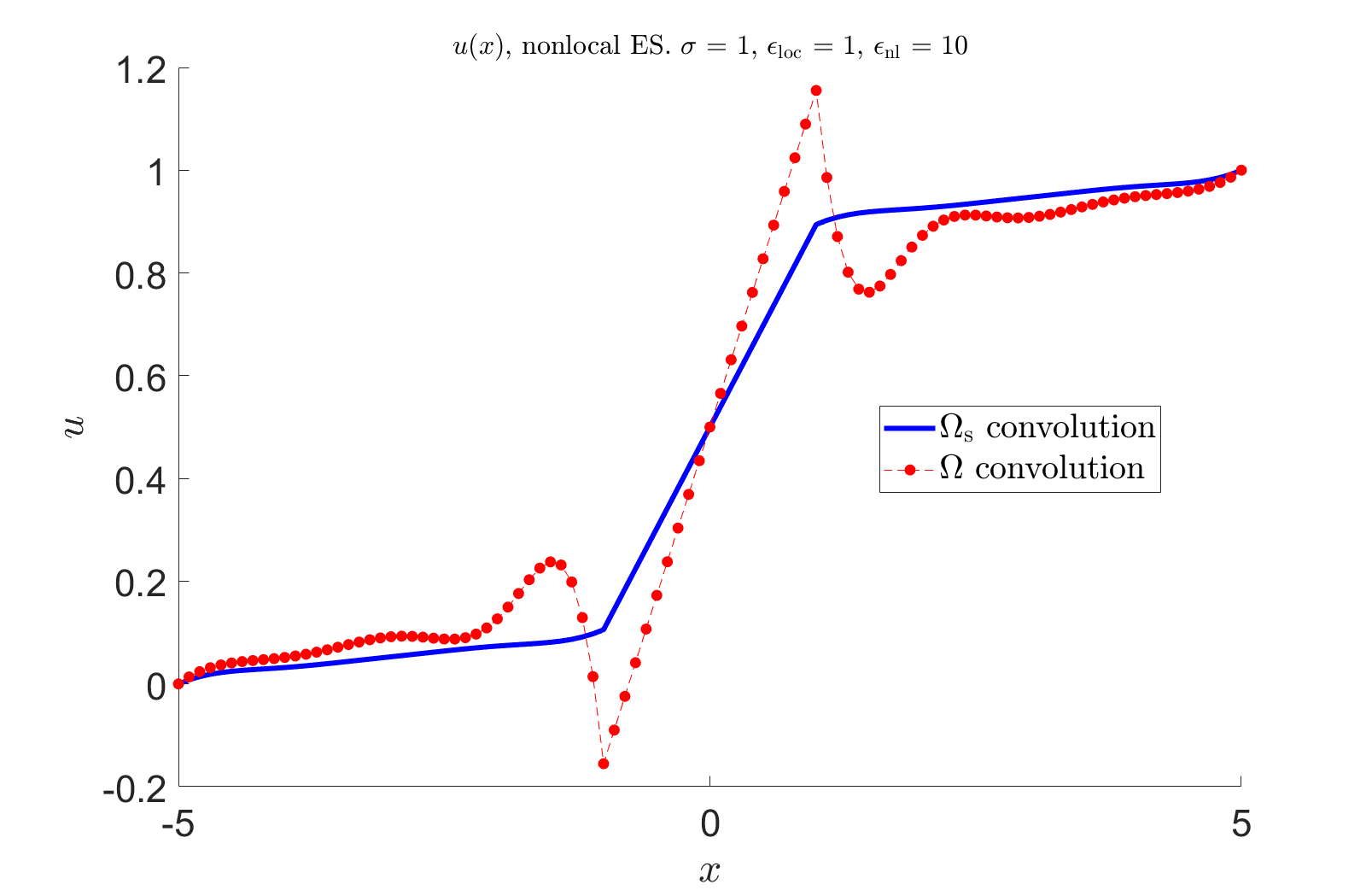}
	\includegraphics[width=0.48\linewidth]{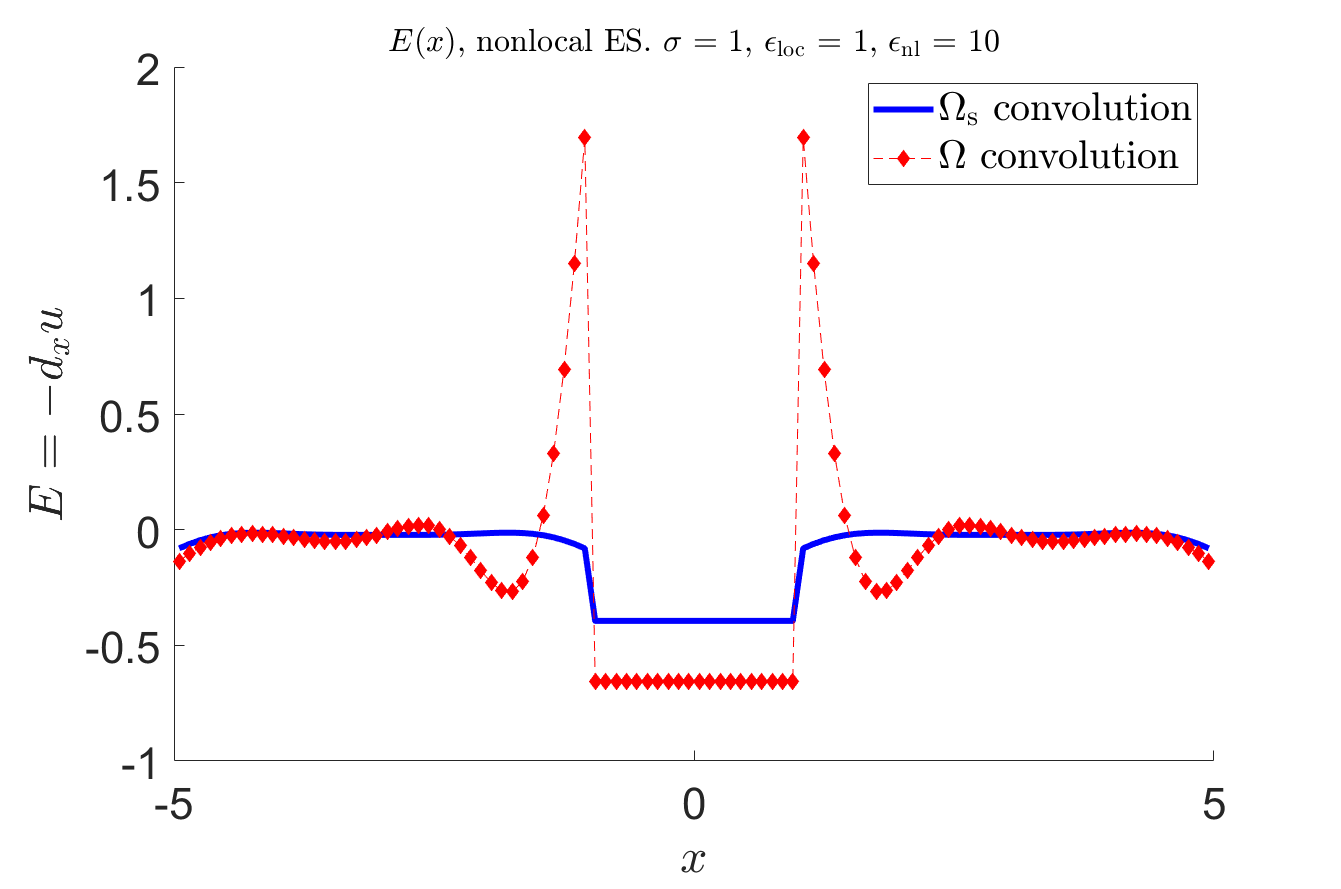}
	\caption{The electrostatic potential and field for a simple
		nonlocal problem in 1D.}
	\label{fig:u-E-nonlocal-1D}
\end{figure}

Importantly, one may challenge 
the tacit assumption that extension of the convolution domain 
to the solute has only mild physical effects. As a simple illustration,
we consider two closely related 1D problems on the interval 
$\Omegarm = (-5, 5)$, with $\Omeganl = (-5, -1) \cup (1, 5)$,
$\Omegarm_l = (-1,1)$:
\begin{equation}\label{eqn:1D-nonlocal-models}
  	E_{1,2}(x) = -d_x u_{1,2}(x); \quad
  	d_x D = 0  \quad(\Rightarrow ~ D = \mathrm{const})
\end{equation}
\begin{equation}\label{eqn:1D-nonlocal-models-bc}
  	D_{1,2}(x_0+) = D_{1,2}(x_0-), ~~ x_0 = \pm 1;
  	\quad
  	u_{1,2}(-5) = 0; ~~ u_{1,2}(5) = 1
\end{equation}
\begin{equation}\label{eqn:1D-D-vs_E-models}
  	D_1(x) = \calE *_{\mathrm{\Omega}} E_1 \equiv
    \int_0^1 \calE(x-x') E_1(x') dx';
    \quad \quad
	D_2(x) = \calE *_{\mathrm{\Omega}_s} E_2 \equiv
	\int_{\mathrm{\Omega}_s} \calE(x-x') E_2(x') dx'
\end{equation}
\begin{equation}\label{eqn:calE-Gaussian-1D}
   \calE(x) \,=\, \exp \left( -\frac{x^2}{2 \sigma^2} \right),
   \quad \quad \sigma = 1
\end{equation}
The two problems differ only in the domain of convolution with
a Gaussian (for convenience) kernel. The results, however, are
significantly different, as evident by the potential and field
plots in \figref{fig:u-E-nonlocal-1D}, and by the fact that
$D_1 = D_1(x) \approx -0.3944$, while $D_2$ is almost double that value,
$\approx -0.6553$.

One concludes that the treatment of convolution, with the
integration domain confined to the region of nonlocality $\Omeganl$,
deserves closer attention. 

From the physical perspective, one caveat is that the dielectric behavior
of solvents in the vicinity of interfaces is not known well enough.
Hence an objection may be raised that the mathematical model need not
be more accurate than the available physical model. 
A possible response to this objection is three-fold. 
First, as the above 1D example illustrates, the treatment of 
nonlocal conditions might affect the solution not mildly and quantitatively, 
but qualitatively; this difference may prove to be critical for
solvation problems, protein folding, and many other problems in
biophysics. Second, it is of computational/mathematical interest 
to solve nonlocal problems efficiently and accurately.
Third, if such accurate solutions are available, one may hope
to start solving inverse problems -- finding the best approximations
to the dielectric function of the solvent near boundaries,
although this ambitious goal is well beyond the scope of this paper.

Similar considerations apply to more complex models with integration kernels
lacking translational invariance near interfaces and thus depending on two
position vectors $\bfr$, $\bfr'$ rather than just on their difference.
%
%%%%%%%%%%%%%%%%%%%%%%%%%%%%%%%%%%%%%
\subsection{The Hildebrandt Localization}
\label{sec:Hildebrandt-Localization}
%%%%%%%%%%%%%%%%%%%%%%%%%%%%%%%%%%%%%
%
An ingenious idea proposed by Hildebrandt \etal 
\cite{Hildebrandt04,Weggler-Rutka-Hildebrandt-JCompPhys10} in the early 2000s
and later enhanced and efficiently implemented by 
by Xie \etal \cite{Xie-SIAM12,Xie-SIAM13,Ying-Xie18}, allows one to convert
nonlocal electrostatic problems to coupled local ones.
This conversion is valid under the simplification assumptions noted below.

The essence of Hildebrandt's approach can be explained as follows
(technical details are available in the literature cited above).
In the nonlocal region $\Omeganl$, we have
\begin{equation}\label{D-eq-calE-conv-grad-u}
    \bfD(\bfr) \,=\, \calE * \bfE \,=\, -\calE * \nabla u   
    ~ \overset{!}{=} ~ -\nabla (\calE * u) \,\equiv\, -\nabla U,
    \quad \quad
    U(\bfr) \overset{\mathrm{def}}{=} \calE * u,
    \quad \bfr \in \Omeganl
\end{equation}
%
%where (importantly) the convolution can be taken over the whole space
%because, by assumption, for all points $\bfr \in \Omegarm_{\infty}$
%the bulk relationship holds.
%
Thus $U$ is the scalar potential for $\bfD$ in $\Omeganl$.
The governing equation for $U$ in $\Omeganl$ therefore is
\begin{equation}\label{eqn:del2-Phi-eq-0}
	-\nabla^2 U(\bfr) \,=\, \gamma \rho, ~~~~~ \bfr \in \Omeganl
\end{equation}
The nonlocality has not yet been eliminated due to the nonlocal
relation between $u$ and $U$.
%(to repeat, these are potentials of the $\bfE$ and $\bfD$ fields, 
%respectively). 
%In Hildebrandt's theory, 
A critical simplification occurs for a specific but physically important
class of convolution kernels. Namely, suppose that there exists
a ``magic'' differential operator $\calL$ whose Green's function
is the given convolution kernel, i.e.
\begin{equation}\label{eqn:calL-calE-eq-delta}
	\calL \calE(\bfr) \,=\, \delta(\bfr)
\end{equation}
Under that assumption, we obtain
\begin{equation}\label{eqn:calL-U-eq-u}
    \calL U \,=\, \calL (\calE * u) ~ \overset{!}{=} ~ 
    (\calL \calE) * u \,=\, \delta * u = u
~~~ \mathrm{in}~\Omeganl,
\end{equation}
where we used the fact that $\calL$ is a differential operator and that 
differentiation can be performed under the convolution sign.
Thus the original nonlocal problem has been reduced to two local
ones, \eqref{eqn:del2-Phi-eq-0}, \eqref{eqn:calL-U-eq-u},
coupled with the local problem within $\Omegaloc$ via the proper 
interface conditions between $\Omeganl$ and $\Omegaloc$.
(Potential splitting noted above also needs to be incorporated in
the actual procedure.)

A key issue is the transformations marked with the exclamation signs
in \eqref{D-eq-calE-conv-grad-u} and \eqref{eqn:calL-U-eq-u}. 
These transformations rely on the property of standard convolution 
over the whole space, but are invalid
when integration is restricted to $\Omeganl$ 
(\sectref{sec:Domain-of-Convolution} and Remark therein).
%Our method does not have this limitation.
%
%%%%%%%%%%%%%%%%%%%%%%%%%%%%%%%%%%%%%%%%%%%%%%%%%%%%%%%%%%%%%%%%%%%%%%%%%%
\section{Construction of Trefftz Functions for Nonlocal Electrostatics}
\label{sec:Trefftz-funcs-nonlocal}
%%%%%%%%%%%%%%%%%%%%%%%%%%%%%%%%%%%%%%%%%%%%%%%%%%%%%%%%%%%%%%%%%%%%%%%%%%
%
%%%%%%%%%%%%%%%%%%%%%%%%%%%%%%%%%%%%%%%%%%%%%%%%%%%%%%%%%%%%%%%%%%%
\subsection{Trefftz Approximations for the Laplace Equation: Harmonic 
Polynomials}\label{sec:Harmonic-polynomials}
%%%%%%%%%%%%%%%%%%%%%%%%%%%%%%%%%%%%%%%%%%%%%%%%%%%%%%%%%%%%%%%%%%%
%
Harmonic polynomials are known to provide an excellent (in some sense, even
optimal \cite{Babuska97}) approximation of harmonic functions 
\cite{Andrievskii87,Babuska97,Bergman66,Melenk99}. The following result is 
cited in \cite{Babuska97}:

\begin{theorem}
	(Szeg\"{o}). Let $\Omegarm \subset \mathbb{R}^2$ be a simply
	connected bounded Lipschitz domain. Let $\tilde{\Omegarm}
	\supset\supset \Omegarm$ and assume that $u \in
	L^2(\tilde{\Omegarm})$ is harmonic on $\tilde{\Omegarm}$. Then there
	is a sequence $\left( u_{p} \right) _{p=0}^{\infty } $ of harmonic
	polynomials of degree $p$ such that
	%%%%%%%%%%
	$$
	\left\| u-u_{p} \right\| _{L^{\infty } (\Omegarm )} ~\le~ c \,
	\exp(-\gamma p) \, \| u \| _{L^{2} (\tilde{\Omegarm})}
	$$
	\begin{equation}
		\| \nabla (u-u_{p}) \| _{L^{\infty } (\Omegarm )} ~\le~ c \,
		\exp(-\gamma p) \, \| u \| _{L^{2} (\tilde{\Omegarm})}
	\end{equation}
	where $\gamma$, $c > 0$ depend only on $\Omegarm$,
	$\tilde{\Omegarm}$.
\end{theorem}

For comparison, the $H^1$-norm error estimate in standard FEM is

\begin{theorem}
	(Ciarlet \& Raviart, Babuska \& Suri
	\cite{Ciarlet72}, \cite{Ciarlet80}, \cite{Babuska94b}). For a family
	of quasiuniform meshes with elements of order \textit{p} and maximum
	diameter $h$, the approximation error in the corresponding finite 
	element space $V^n$ is
	$$
	\inf_{v \in V^n} \| u-v \| _{H^1 (\Omegarm )} ~=~ C h^{\mu-1}
	p^{-(k-1)} \, \| u \| _{H^k (\Omegarm )}
	$$
	where $\mu = \min(p+1, k)$ and $c$ is a constant independent of $h$,
	$p$, and $u$.
\end{theorem}

Thus, for a fixed polynomial order $p$, the FEM and harmonic approximation 
errors are similar \cite{Babuska97}; however, the FEM approximation 
is realized in a much wider space containing \textit{all}
polynomials up to order $p$, not just harmonic ones. For solving
the Laplace equation, the standard FE basis set can thus be viewed 
as having substantial redundancy that is eliminated by using the
harmonic basis.

\textit{These results and observations motivate our development of Trefftz
functions in the nonlocal case; we conjecture that similar 
excellent approximation properties will hold.}
This is supported by the numerical experiments of 
\sectref{sec:Trefftz-functions-nonlocal-2D}.
%
%%%%%%%%%%%%%%%%%%%%%%%%%%%%%%%%%%%%%%%%%%%%%%%%%%%%%%%%%%%%%%%%%%%
\subsection{Trefftz Functions for Nonlocal 
Problems}\label{sec:Trefftz-functions-nonlocal-setup}
%%%%%%%%%%%%%%%%%%%%%%%%%%%%%%%%%%%%%%%%%%%%%%%%%%%%%%%%%%%%%%%%%%%
%
%%%%%%%%%%%%%%%%%%%%%%%%%%%%%%%%%%%%%%%%%%%%%%%%%%%%%%%%%%%%%%%%%%%
\subsubsection{1D Construction: Fixing 
Ideas}\label{sec:Trefftz-functions-nonlocal-1D}
\begin{figure}
	\centering
	\includegraphics[height=2.5cm]{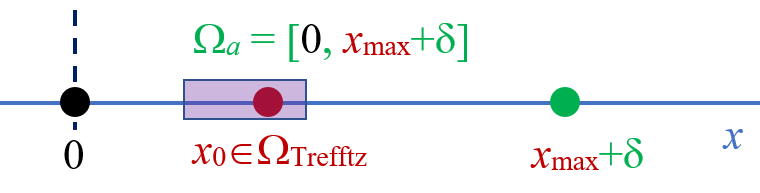}
	\caption{Construction of Trefftz functions in a (small) 1D
		domain $\Omegarm_{\mathrm{Trefftz}} \equiv (x_{\min}, x_{\max})$ 
		indicated with a shaded rectangle. 
		Only the $E$ field in $[0, x_{\max}+\delta]$ contributes
		to $D(x)$ in $\Omegarm_{\mathrm{Trefftz}}$. Potential
		$u(x)$ is approximated by a set of basis functions (e.g. Taylor 
		polynomials), and the expansion coefficients are determined
		by imposing, to a certain order, zero divergence conditions for $D$ 
		in $\Omegarm_{\mathrm{Trefftz}}$. 
		See text for further details.}
	\label{fig:setup-u-Trefftz-1D}
\end{figure}

%%%%%%%%%%%%%%%%%%%%%%%%%%%%%%%%%%%%%%%%%%%%%%%%%%%%%%%%%%%%%%%%%%%
%
For simplicity of exposition, assume that the support 
of the convolution kernel is finite: 
diam(supp~$\calE$)~$= \delta$ % \ll$~diam($\Omegarm$).
(if needed, all physically meaningful kernels in nonlocal electrostatics
can be truncated with an exponentially small error).
To fix ideas, let us start with a 1D setup (\figref{fig:setup-u-Trefftz-1D}).
Let $x > 0$ represent a nonlocal medium, and $x = 0$ be either the
boundary of the computational domain or an interface boundary between
local ($x < 0$) and nonlocal ($x > 0$) media. 

\textit{The objective is to derive Trefftz functions} $u_{\mathrm{Trefftz}}(x)$ 
	for a (small) interval
	$\Omegarm_{\mathrm{Trefftz}} = (x_{\min}, x_{\max})$, indicated
in \figref{fig:setup-u-Trefftz-1D} with a shaded purple rectangle.
Assume that $x_{\min} < \delta$
(the case $x_{\min} > \delta$ is handled analogously and is slightly simpler).
One observes that the $D$ field within $\Omegarm_{\mathrm{Trefftz}}$
may depend only on the $E$ field within $[0, x_{\max}+\delta]$.

With this in mind, introduce a basis set $\{u_{\alpha}\}_{\alpha=0}^n$,
$n$ being an adjustable parameter. An obvious, but not only, choice
is the Taylor polynomials $u_{\alpha}(x) = (x-x_0)^{\alpha}$, where
$x_0$ is an arbitrary point in $\Omegarm_{\mathrm{Trefftz}}$.
Then we have the respective functions
\begin{equation}\label{eqn:ED-basis}
	E_{\alpha}(x) = -d_x u_{\alpha}(x);  \quad
	D_{\alpha}(x) = \int_0^{x_{\max}+\delta} \calE(x-x') E_{\alpha}(x') \, dx'
\end{equation}
Looking for $u_{\mathrm{Trefftz}}$ as a linear combination
%
%\begin{equation}\label{eqn:u-trefftz-sum-c-u}
$
	u_{\mathrm{Trefftz}}(x) \,=\, \sum\nolimits_{\alpha=0}^n c_{\alpha} 
	u_{\alpha}(x)
$
%\end{equation}
%
and imposing conditions for the ``1D divergence'' of $D$
\begin{equation}\label{eqn:dx-beta-eq-0}
	d_x^{\beta} \, \mathrm{div} \, D_{\mathrm{Trefftz}}(x_0) \,\equiv\, 
	\sum_{\alpha=0}^n c_{\alpha} \, d_x^{\beta} D_{\alpha}(x_0) = 0,
	\quad \beta = 1,2, \ldots m
\end{equation}
(where $d_x^{\beta}$ is shorthand for $d^{\beta}/dx^{\beta}$),
one obtains Trefftz solutions whose number depends on the adjustable 
parameters $n, m$. I call these solutions
\textit{pseudoharmonic functions} -- by analogy with 
harmonic polynomials and ``generalized harmonic polynomials'' 
\cite{Melenk99,Hiptmair-Trefftz-survey16}.

\begin{figure}
	\centering
	\includegraphics[height=1.45in]{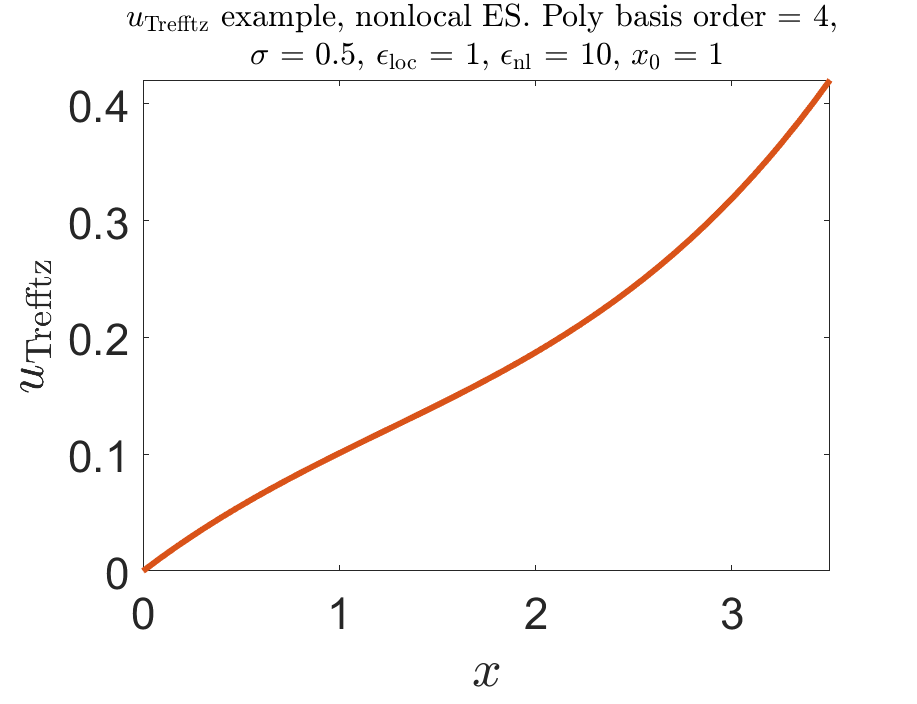}
	\includegraphics[height=1.45in]{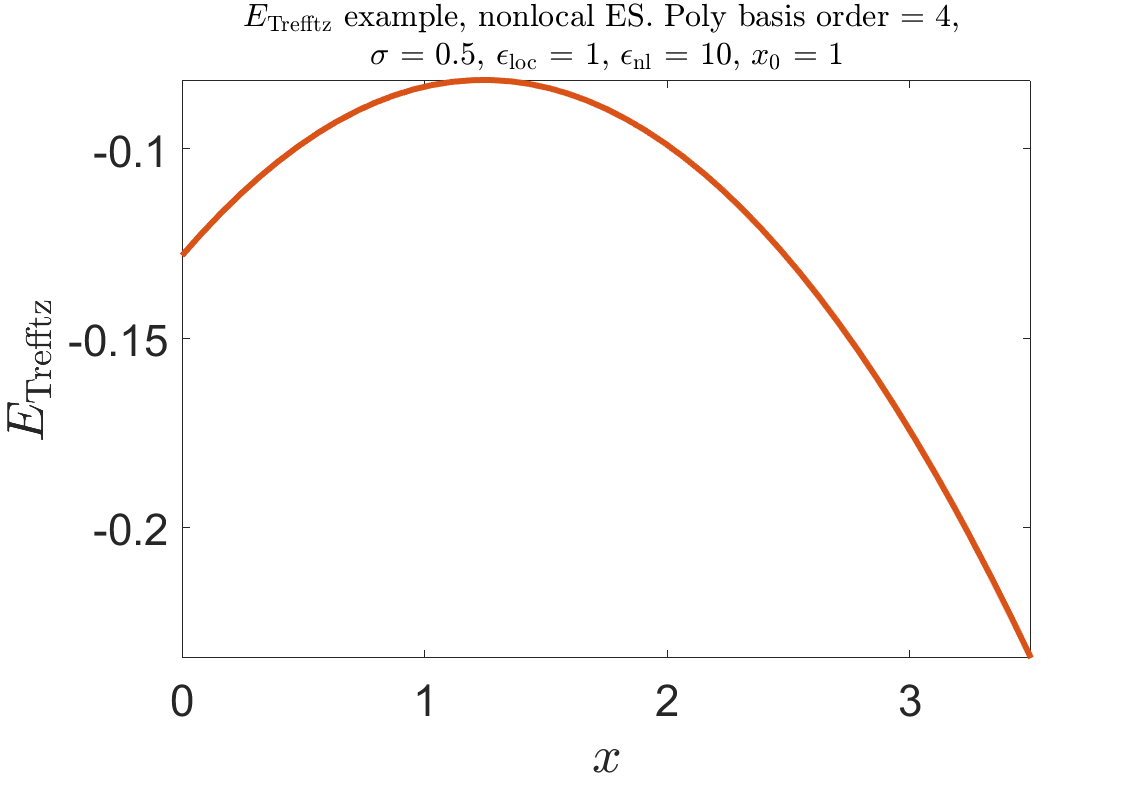}
	\includegraphics[height=1.45in]{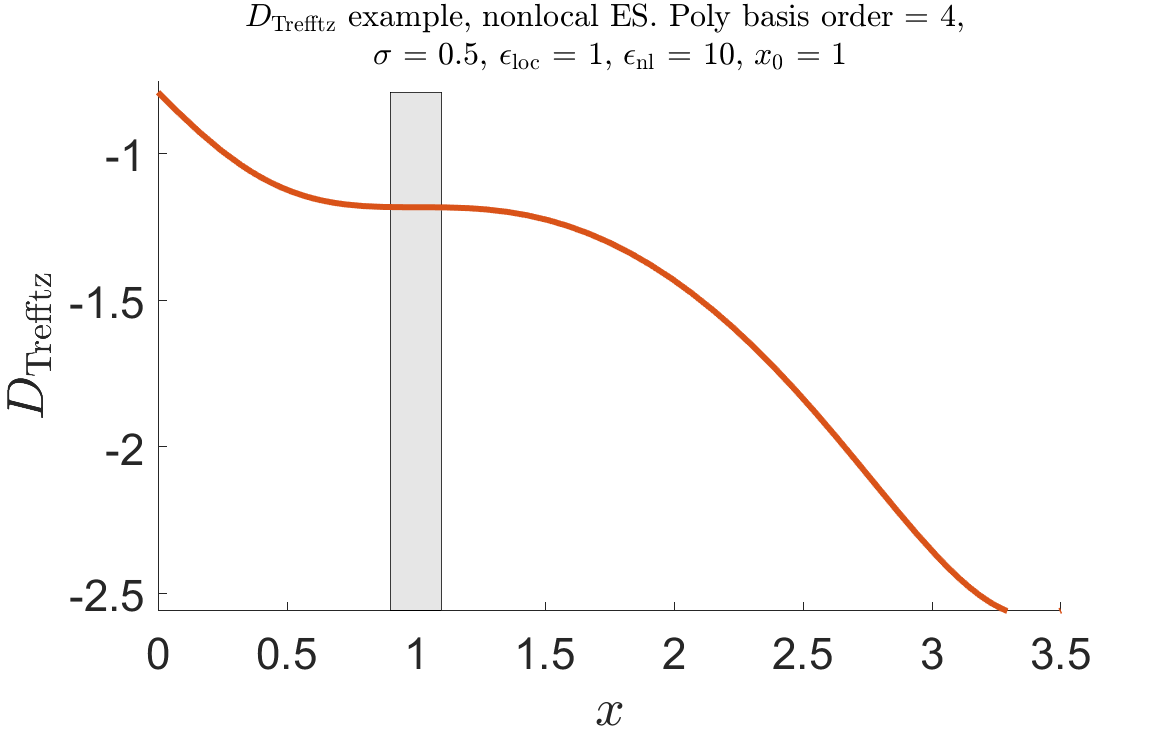}
	\caption{
		One	of the three Trefftz functions $u(x)$ (left),
		with the corresponding $E(x)$ (middle) and $D(x)$ (right).
		Note the flatness of $D(x)$ in the shaded vicinity of
		$x_0 = 1$. Parameters: $\epsloc = 1$, $\epsnl = 10$,
		$x_0 = 1$, $n=4$, $m=1$; the Gaussian kernel with $\sigma = 0.5$.}
	\label{fig:uED-Trefftz-example-1D}
\end{figure}

Plotted in \figref{fig:uED-Trefftz-example-1D} is one
of the three pseudoharmonic functions for $\epsloc = 1$, $\epsnl = 10$,
$x_0 = 1$, $n = 4$, $m=1$, and the Gaussian kernel 
$\calE(x) = \exp(-x^2/(2 \sigma^2))$, $\sigma = 0.5$. 
Note the flatness of the $D$ curve around $x_0 = 1$ (shaded area),
as expected. Also note that in the nonlocal domain $D$ is not proportional 
to $E$.

This construction extends in a natural way to higher dimensions;
the technical difficulties of computing convolution integrals
exactly or approximately are overcome with the help
of symbolic algebra (see the 2D example below).

While Trefftz sets are interesting in their own right,
their main application is in the Trefftz-DG or Trefftz-FLAME context:
then, $\Omegarm_{\mathrm{Trefftz}}$ will contain a given finite element, 
or a grid ``molecule''. The Trefftz sets will in general vary from element to 
element
or stencil to stencil, but that is a natural feature of both DG and FLAME.
Importantly, \textit{one can then ``forget'' about nonlocality:
information about it is built into the Trefftz functions,
which can be used the same way as if the problem were local}.

At an interface between nonlocal and local media, one additional
step is required. Pseudoharmonic functions in the nonlocal medium
must be glued with their local counterparts (typically,
harmonic polynomials) via the matching conditions. 

The pseudoharmonic functions constructed here should,
strictly speaking, be classified as 
``quasi-Trefftz''~\cite{Imbert-Gerard-arxiv19,Imbert-Gerard-arxiv20a,
	Imbert-Gerard-arxiv20b}, since they satisfy
the underlying equation to a desired order rather than exactly.
I still retain the plain term ``Trefftz functions'' for brevity,
at the expense of a mild abuse of the terminology.\footnote{One may conjecture 
that the 
errors inherent in quasi-Trefftz functions can be ignored if they are of higher 
order than other numerical errors in a particular method.
The ``quasi'' prefix can be restored in cases where it would make a difference.}
%
%%%%%%%%%%%%%%%%%%%%%%%%%%%%%%%%%%%%%%%%%%%%%%%%%%%%%%%%%%%%%%%%%%%
\subsubsection{2D Examples}\label{sec:Trefftz-functions-nonlocal-2D}
\begin{figure}
	%	\centering
	\includegraphics[height=1.4in]{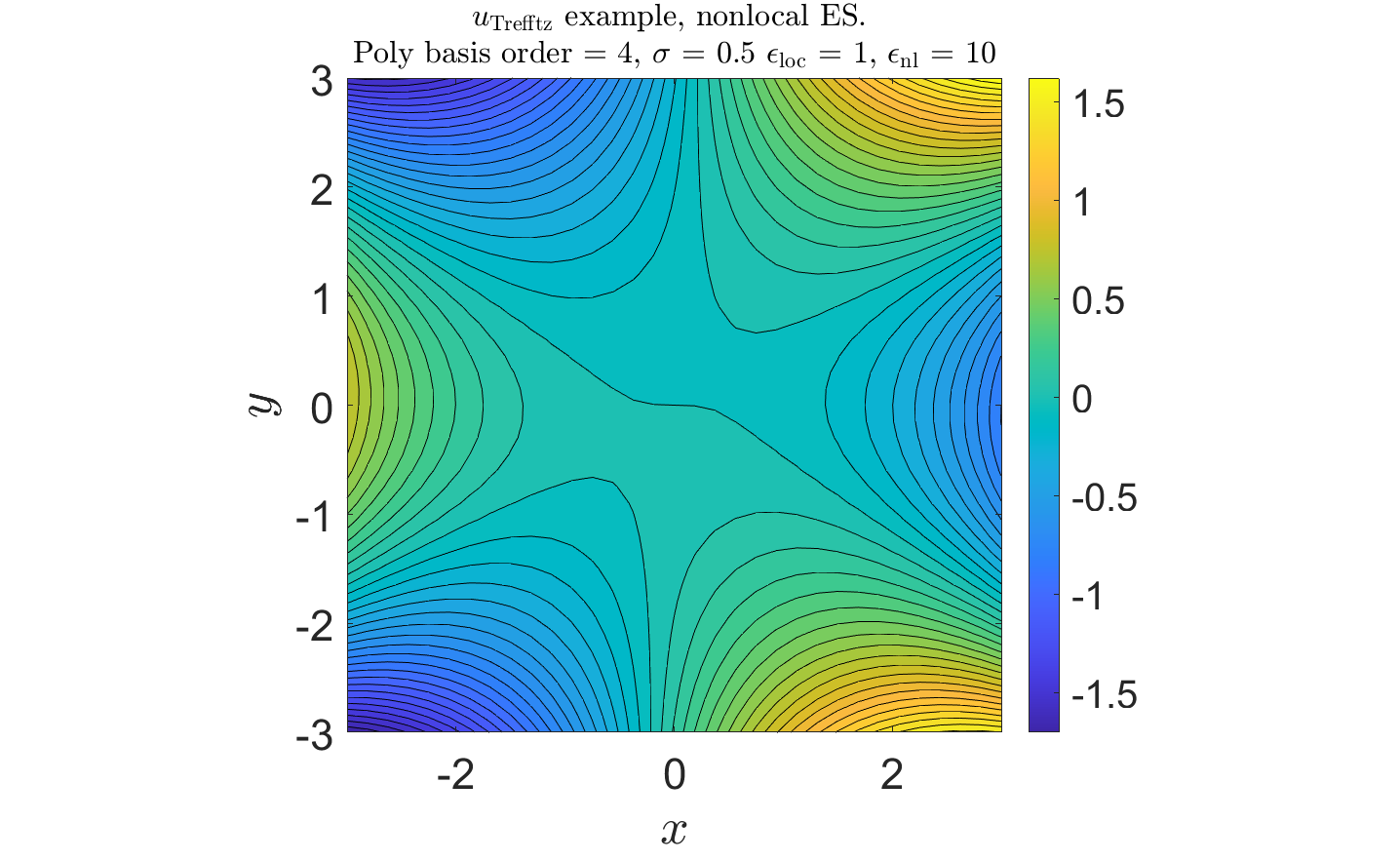}
	\hbox{\hspace{-2ex}}	
	\includegraphics[height=1.4in]{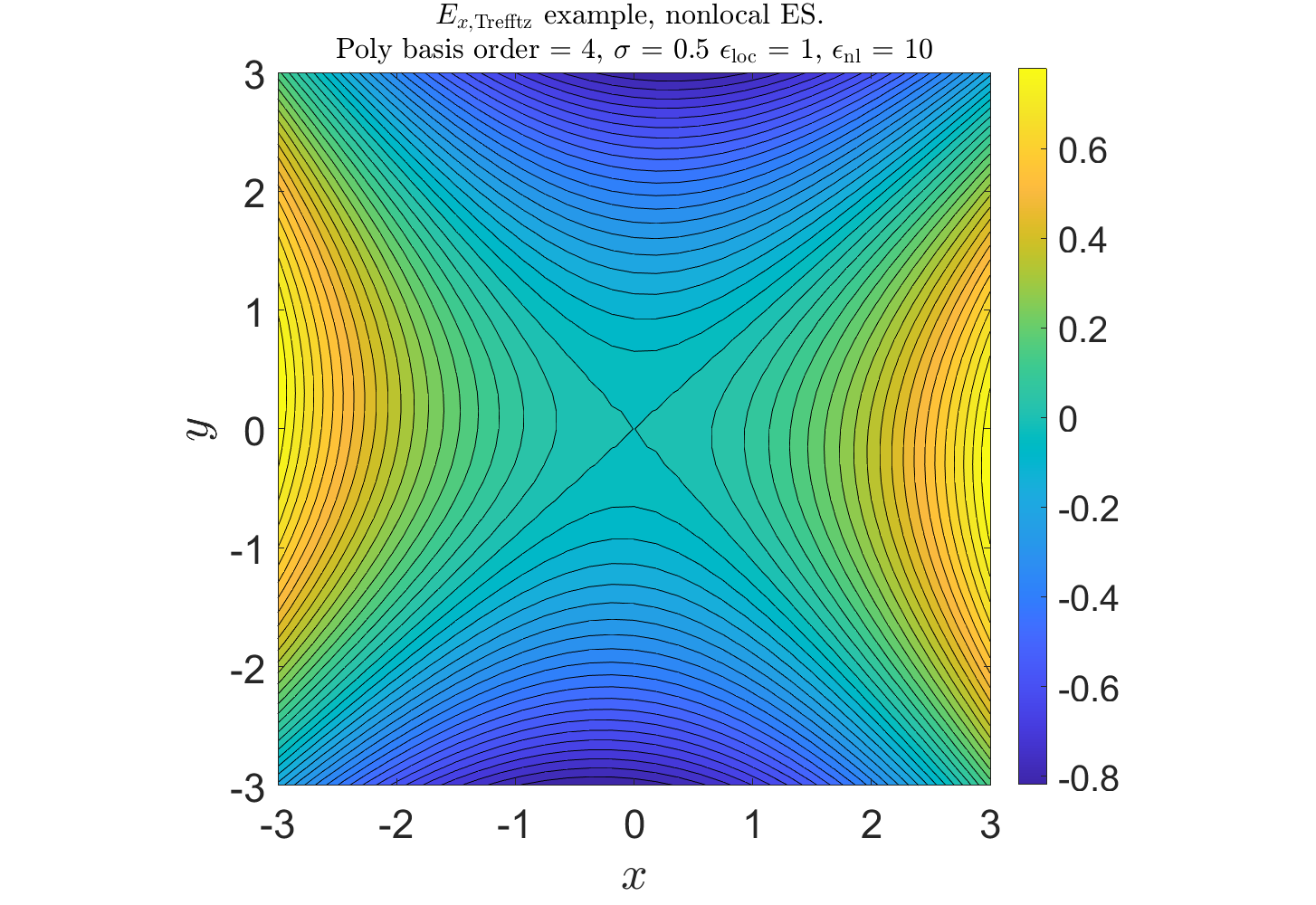}
	\hbox{\hspace{2ex}}	
	\includegraphics[height=1.4in]{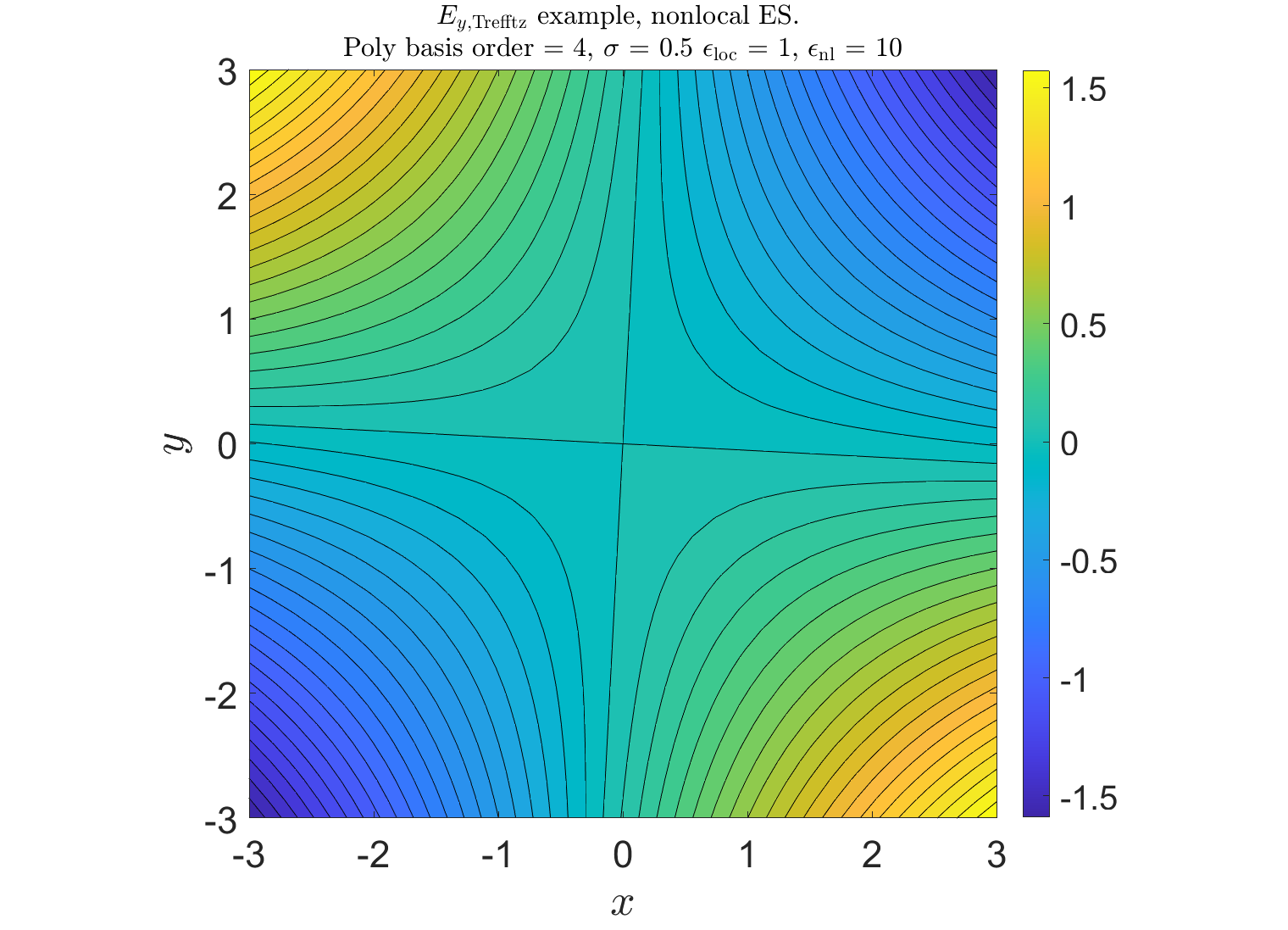}\\
	\hbox{\hspace{1mm}}
    \includegraphics[height=1.4in]{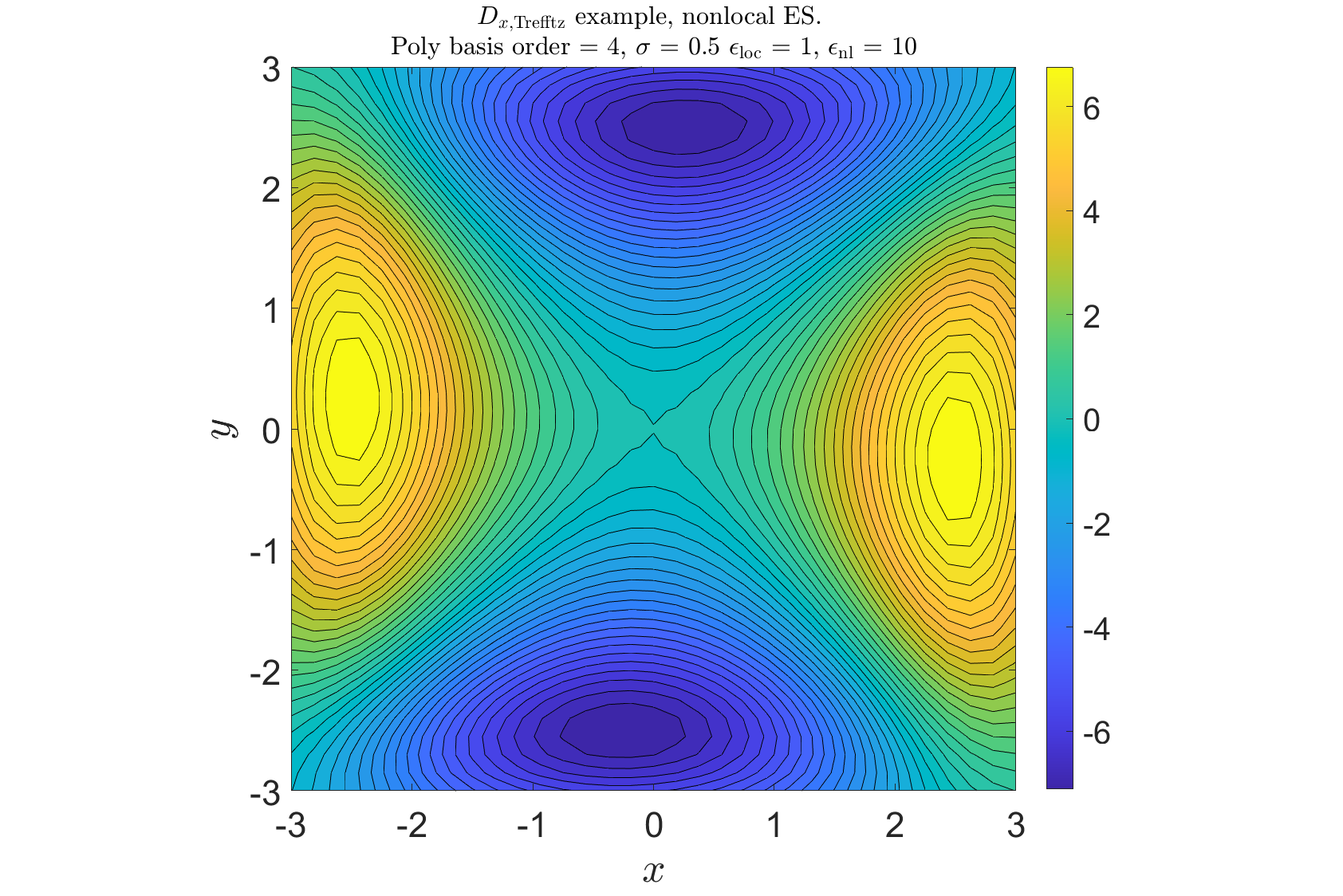}
    \includegraphics[height=1.4in]{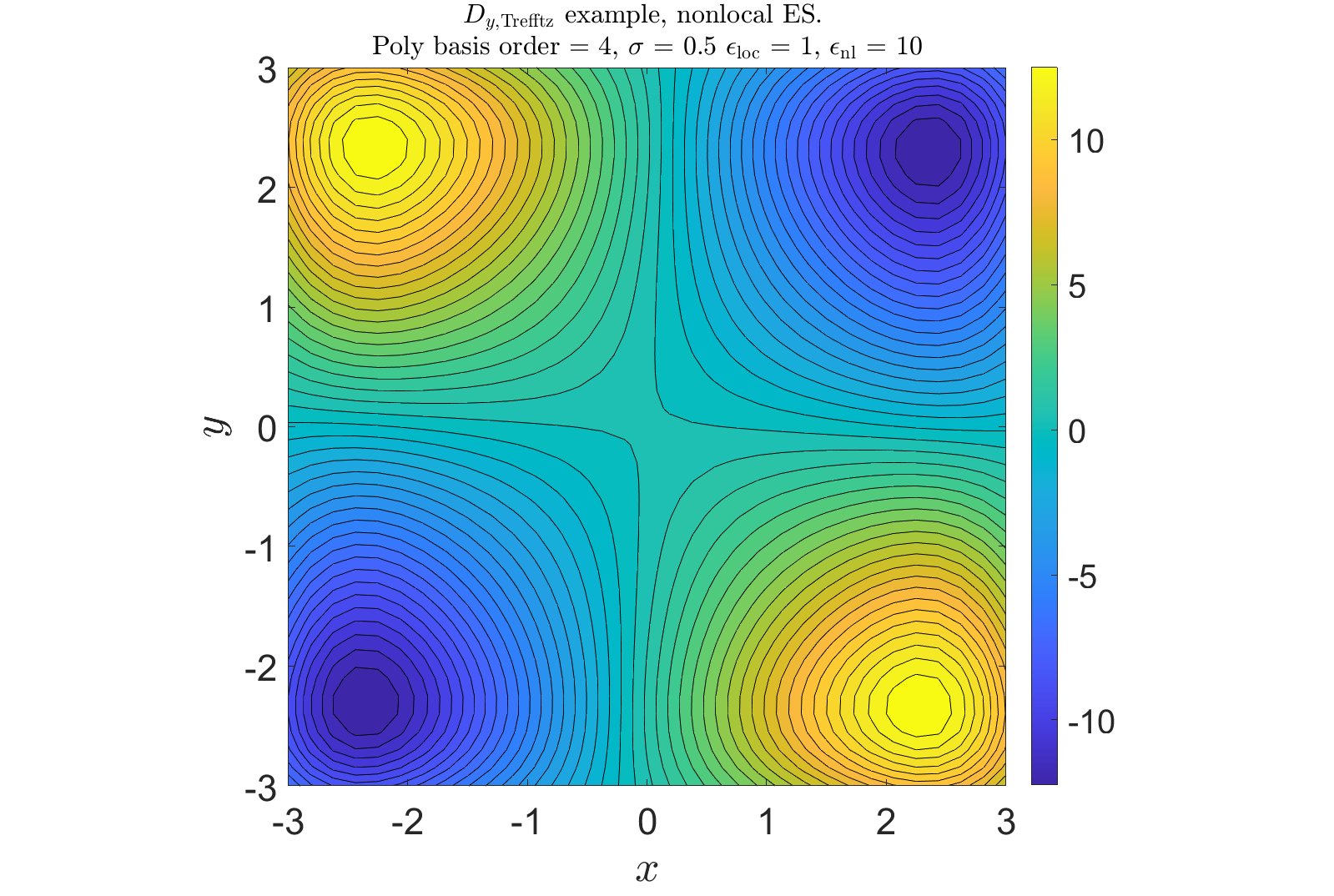}
	\includegraphics[height=1.4in]{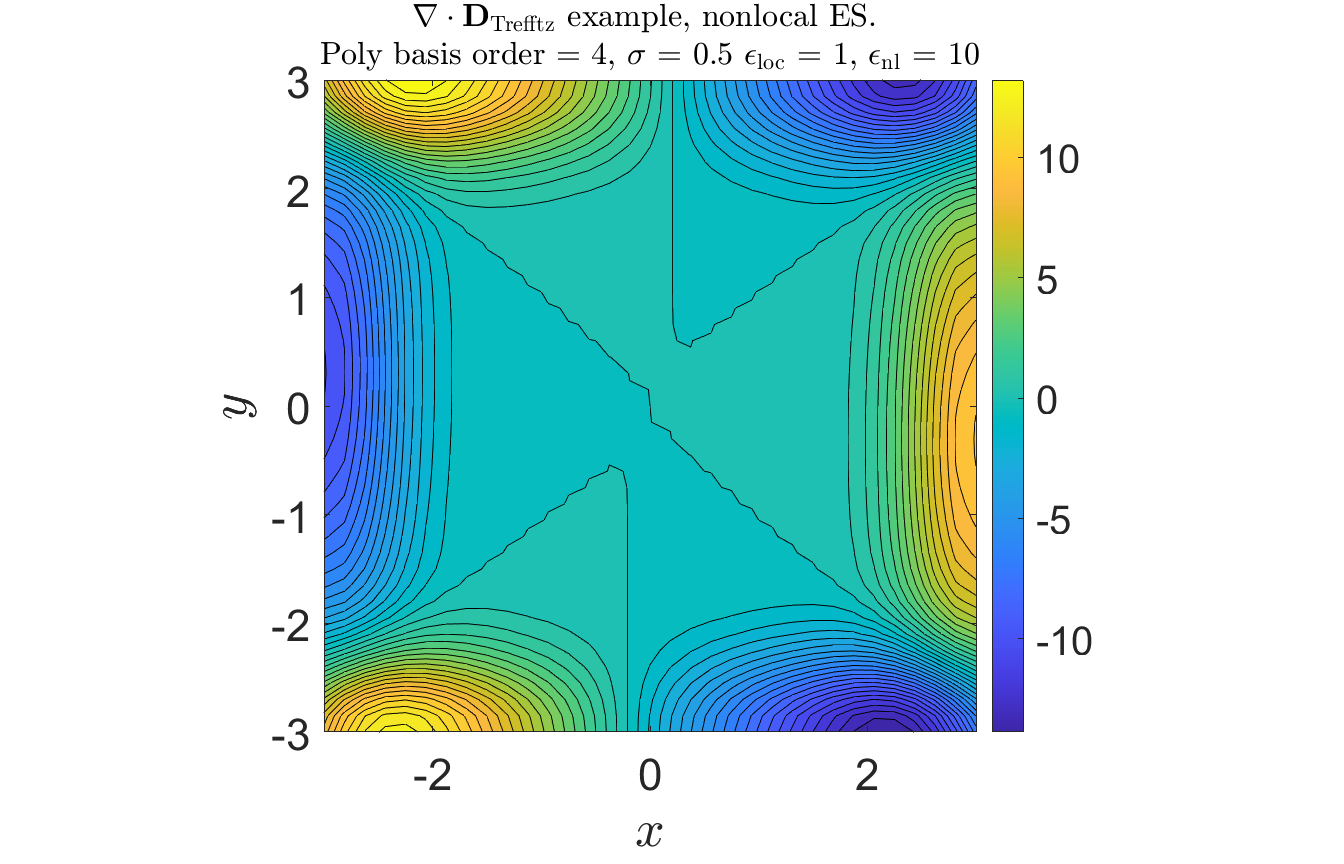}
%	\vspace{-10mm}
	\caption{A Trefftz function ``in the bulk'' of a nonlocal medium.
		Polynomial order for the potential $u(x, y)$ is $n_{\max} = 4$.
		Gaussian kernel (in $x, y$) with $\sigma = \frac12$;
		$\epsloc = 1$, $\epsnl = 10$. Zero-divergence condition
		for $\bfD$ imposed at $(x_0, y_0) = (0,0)$ with partial derivatives 
		up to order	$m = 2$. Top row: potential $u_{\mathrm{Trefftz}}$
	and the corresponding $E_x$, $E_y$. Bottom row: $D_x$, $D_y$,
     and $\nabla \cdot \bfD$.}
	\label{fig:u-E-D-Trefftz-bulk}
\end{figure}

%%%%%%%%%%%%%%%%%%%%%%%%%%%%%%%%%%%%%%%%%%%%%%%%%%%%%%%%%%%%%%%%%%%
%
In 2D, construction of pseudoharmonic functions is conceptually the same
as in 1D. First, consider approximations ``in the bulk'' (i.e. away
from any interfaces). As an illustrative example,
let the domain $\Omegarm_a$, analogous to the one
in \figref{fig:setup-u-Trefftz-1D}, be $[-6 \sigma, 6 \sigma]^2$,
with $\sigma = \frac12$. Let the electrostatic potential $u$ be
approximated in $\Omegarm_a$ by the set of all polynomials
in $x, y$ up to order $n_{\max} = 4$, leading to the respective
fields $\bfE = -\nabla u$ and $\bfD = \calE *_{\Omega_a} \bfE$,
where $\calE$ is a Gaussian kernel in $x, y$ with the parameter $\sigma$;
$\epsloc = 1$, $\epsnl = 10$. Then impose conditions analogous to
\eqref{eqn:dx-beta-eq-0}, with partial derivatives taken up to order
$m = 2$ at $(x_0, y_0) = (0,0)$.
One of the Trefftz functions generated this way is shown in
\figref{fig:u-E-D-Trefftz-bulk}.

\begin{figure}
	\centering
	\hbox{\hspace{3ex}}
	\includegraphics[width=0.75\linewidth]{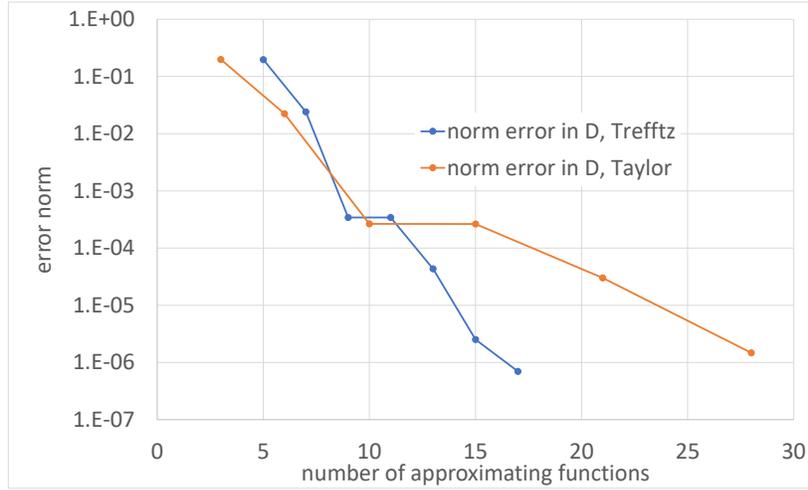}
	\vspace{-10mm}
	\caption{Error norms vs. the number	of approximating functions. 
		Trefftz approximation of $\bfD_{\mathrm{test}}$ is seen to
	converge much more rapidly than the standard Taylor expansion.
	$D_{\mathrm{test}, x} = \partial_y v_{\mathrm{test}}$; 	
	$D_{\mathrm{test}, y} = -\partial_x v_{\mathrm{test}}$;
	$v_{\mathrm{test}} = \sin x \exp y + \exp(-x) \cos y$.
}
	\label{fig:error-norm-vs-num-funcs}
\end{figure}

Exponential convergence of these Trefftz approximations with respect
to the number of approximating functions (which in turn depends on
the chosen polynomial order $n_{\max}$) is illustrated by 
\figref{fig:error-norm-vs-num-funcs}, in comparison with standard Taylor 
expansions. The approximation error is defined as
\begin{equation}\label{eqn:approx-error}
	\epsilon_a(V) ~\overset{\mathrm{def}}{=}
	\inf_{u_a \in V} \| \bfD_a(u_a) - \bfD_{\mathrm{test}} \|_V
\end{equation}
where
\begin{itemize}
	\item
	$\bfD_a(u_a) = -\calE *_{\Omegarm_a} \nabla u_a$.
	\item $\Omega_t$ is the ``target'' domain for the condition
	$\nabla \cdot \bfD \approx 0$. In the case of Trefftz
	approximations, $\Omega_t \equiv \Omegarm_{\mathrm{Trefftz}}$.
	In the numerical example, $\Omega_t = (-\sigma, \sigma)^2$.
	\item 
	$V = V(\Omega_t)$ is the approximating space, spanned either by
	the constructed Trefftz set or, alternatively, by Taylor
	polynomials of orders up to $n_{\max}$. In the case of
	Taylor polynomials, the $x$ and $y$ components of $\bfD$
	are approximated separately, which arguably is an unfair advantage;
	yet convergence is still less rapid than in the case of Trefftz 
	approximations.
	\item
	In either case,  $\| \bfD \|_V^2 = \mathrm{area}^{-1}(\Omegarm_t)
	\int_{\Omegarm_t} |\bfD|^2 \, dx \, dy$
	\item A zero-divergence	test function $\bfD_{\mathrm{test}}(x,y)$ 
	is chosen as an example:
		$D_{\mathrm{test}, x} = \partial_y v_{\mathrm{test}}$; 	
		$D_{\mathrm{test}, y} = -\partial_x v_{\mathrm{test}}$;
		$v_{\mathrm{test}} = \sin x \exp y + \exp(-x) \cos y$.
\end{itemize}

% Code: Trefftz_funcs_nonlocal_ES_at_point_sym.m

\begin{figure}
	\centering
%	\hbox{\hspace{3ex}}
	\includegraphics[width=0.33\linewidth]{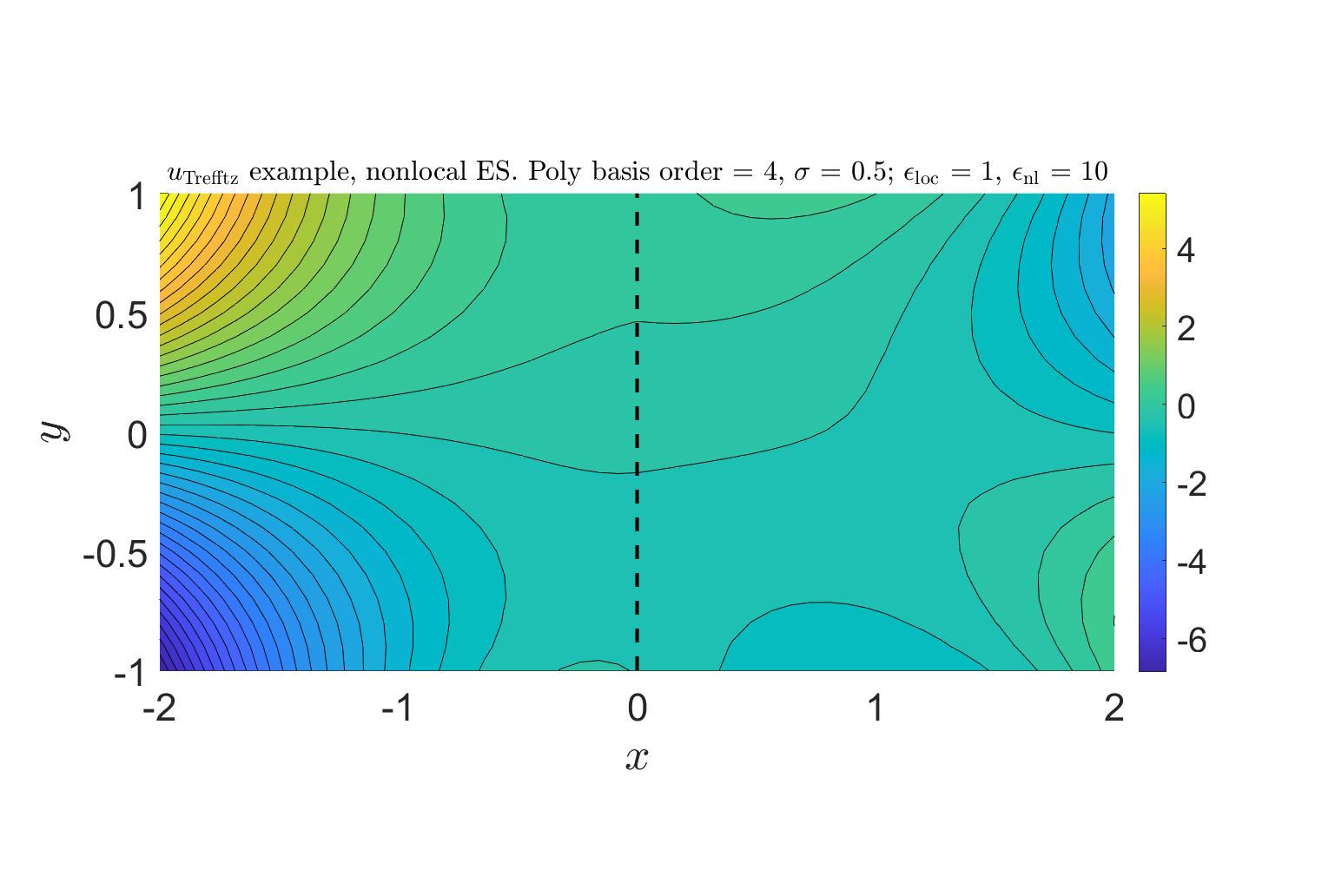}
		\hbox{\hspace{-3ex}}
	\includegraphics[width=0.33\linewidth]{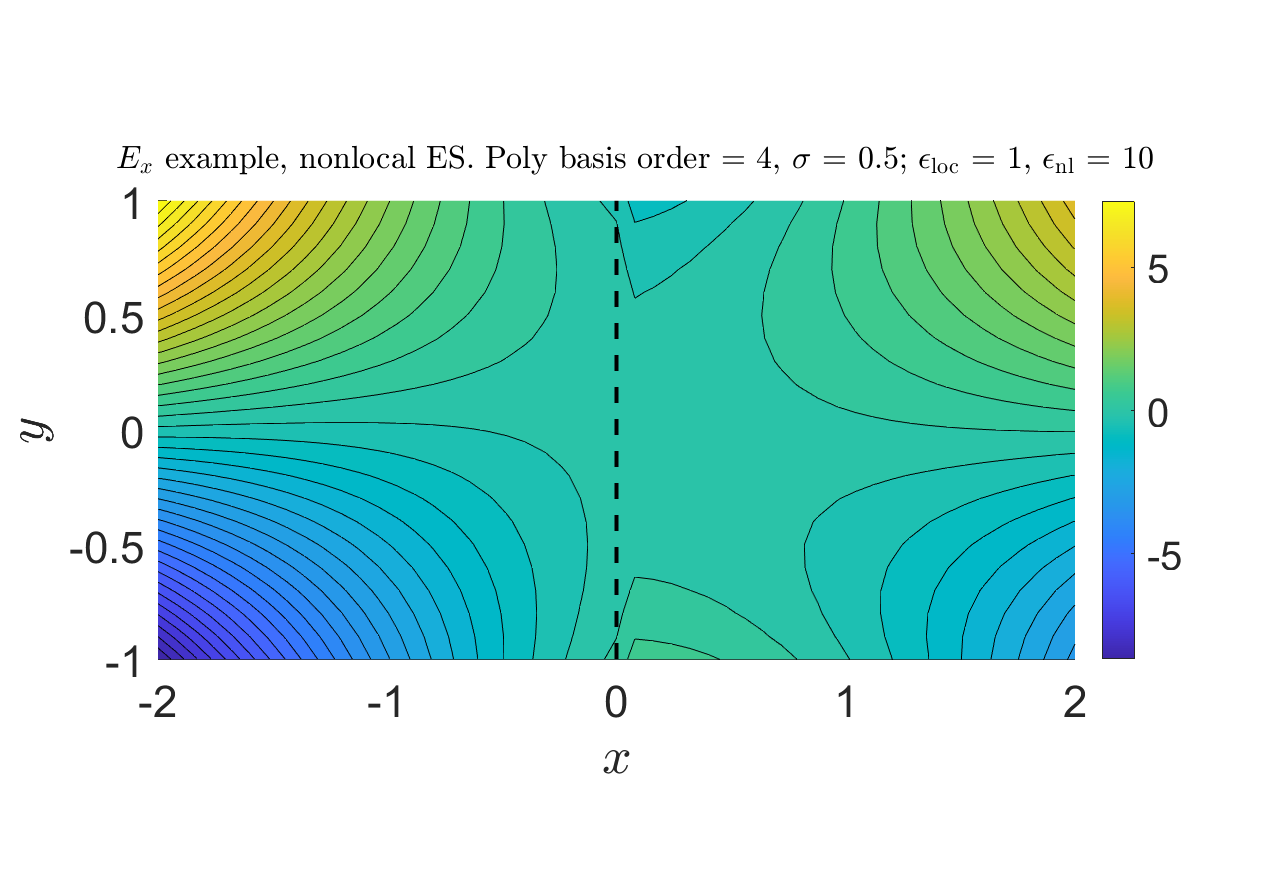}
	\includegraphics[width=0.33\linewidth]{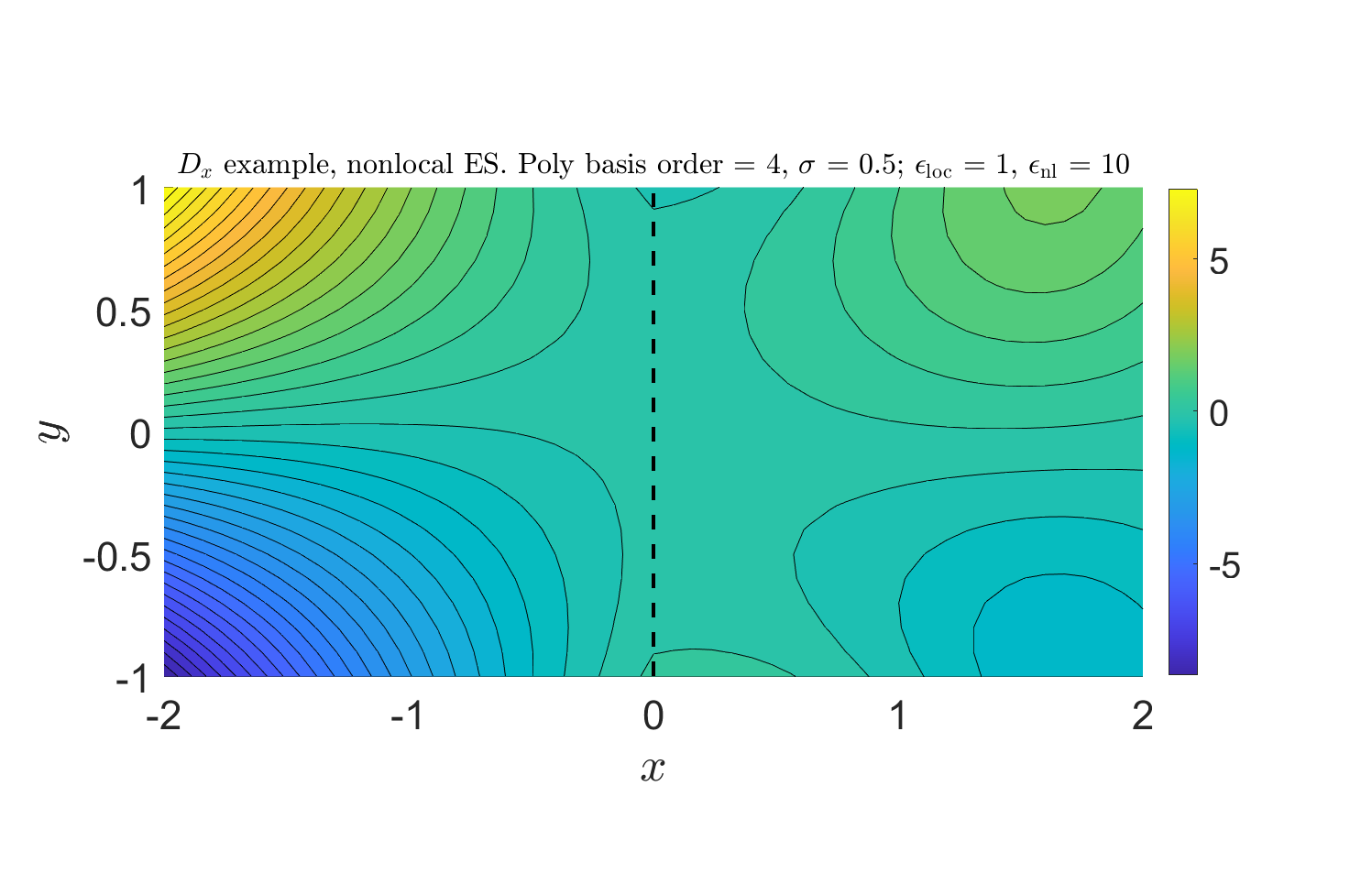}
	\vspace{-10mm}
	\caption{
		Examples of 2D pseudoharmonic functions. Local and nonlocal
		regions: $x < 0$ and $x > 0$, respectively. Taylor order $n = 4$, 
		$p_{\max} = 2$. $\epsloc = 1$, $\epsnl = 10$. Left: $u$, middle: $E_x$,
		right: $D_x$. Note that in the nonlocal domain $\bfD$ is not 
		proportional
		to $\bfE$.}
	\label{fig:u-Trefftz-2D}
\end{figure}

A more complicated example of a pseudoharmonic function in 2D is presented in 
\figref{fig:u-Trefftz-2D}. The kernel is Gaussian, with $\sigma = 0.5$;
$\Omegarm_{\mathrm{Trefftz}}$ is located at the origin,
and the example also features a straight interface boundary $x = 0$
between the local (left) and nonlocal (right) media.
In the local domain, the Trefftz functions are harmonic polynomials
up to order 4. At the interface, the Trefftz functions satisfy
the usual matching conditions for $u$ and $\bfD \cdot \bfn$;
these are cumbersome but straightforward to impose using symbolic
algebra; this only needs to be done once for any given integration kernel.

% Trefftz_funcs_nonlocal_ES_straight_bdry

In practice, the computation of pseudoharmonic functions simplifies greatly if 
the convolution kernel $\calE$ is separable or can be
approximated as a combination of (a small number of) separable functions,
since double integration reduces to a product of single integrals.
Kernels commonly accepted in nonlocal electrostatics
have the form of the Yukawa potential
%
%\begin{equation}\label{eqn:Yukawa-kernel}
$
	\calE(\bfr, \lambda) \,=\,
	\frac{1}{[4 \pi] \lambda^2} \, \exp \left(- \frac{r}{\lambda} \right)
$
%\end{equation}
%
(single-pole in reciprocal space), where the $4 \pi$ factor
in the square brackets is present in the SI system but not in the Gaussian 
system, and parameter $\lambda$ defines the range of nonlocal interactions. 
Such singular kernels are not separable, and one has
several options for finding pseudoharmonic functions:

(i) Consider separable approximations of the kernel; 
there is extensive literature on canonical tensor decomposition
(\cite{Khoromskij10,Biagioni12,Reynolds-Beylkin17} 
and references therein). 

(ii) Precompute the Trefftz basis numerically and accurately;
note that for any given kernel this needs to be done once.
We previously implemented and used numerical Trefftz functions 
\cite[Sect~5.1,~5.2]{Dai-Webb11} --
this was done for local electrostatics, but the principle remains the same.

%%%%%%%%%%%%%%%%%%%%%%%%%%%%%%%%%%%%%%%%%%%%%%%
\section{Conclusion}\label{sec:Conclusion}
%%%%%%%%%%%%%%%%%%%%%%%%%%%%%%%%%%%%%%%%%%%%%%%
%
Trefftz functions, which, by definition, satisfy locally the underlying
differential equation and applicable interface boundary conditions,
tend to provide highly accurate approximations of the solution.
This has led to the development of high-order Trefftz-FD schemes
and Trefftz-DG methods with exponential convergence.
Proposed in this paper is a way to generate Trefftz functions
for nonlocal problems. Examples are presented, and exponential accuracy
is demonstrated numerically. 

In future research, Trefftz approximations are intended to be applied
to problems of nonlocal electrostatics, important in biophysical and 
macromolecular simulation. It is anticipated that these approximations
will help to alleviate the ``curse of nonlocality'' -- much higher
computational cost in comparison with local problems.

\bibliographystyle{plain}
\bibliography{Trefftz_funcs_nonlocal_ES}

\end{document}